\documentclass[pdflatex,sn-mathphys-num]{sn-jnl}

\usepackage{graphicx}%
\usepackage{multirow}%
\usepackage{amsmath,amssymb,amsfonts}%
\usepackage{amsthm}%
\usepackage{mathrsfs}%
\usepackage[title]{appendix}%
\usepackage{xcolor}%
\usepackage{textcomp}%
\usepackage{manyfoot}%
\usepackage{booktabs}%
\usepackage{graphicx}
\usepackage{float}        
\usepackage{subcaption}   
\usepackage{listings}%
\usepackage{enumitem}

\theoremstyle{thmstyleone}%
\newtheorem{theorem}{Theorem}
\theoremstyle{thmstyletwo}%
\newtheorem{remark}[theorem]{Remark}%
\newtheorem{corollary}[theorem]{Corollary}%
\theoremstyle{thmstylethree}%

\begin{document}

\title[On the Exact Distribution]{On the Exact Distribution of the Sum of Two CIR Processes}


\author*[1]{\fnm{Bilgi} \sur{Yilmaz}}\email{bilgiyilmaz07@gmail.com \&bilgiyilmaz@ksu.edu.tr}

\author[2]{\fnm{A.Alper} \sur{Hekimoglu}}\email{ahekimsfe@gmail.com}
\equalcont{These authors contributed equally to this work.}


\affil*[1]{\orgdiv{Department Mathematics/Faculty of Science}, \orgname{Kahramanmaraş Sütçü İmam University}, \orgaddress{\city{Kahramanmaraş}, \postcode{46050}, \country{Türkiye}}}

\affil[2]{\orgdiv{Model Validation Department}, \orgname{European Investment Bank}, \orgaddress{\street{98-100, Boulevard Konrad Adenauer}, \city{Luxembourg}, \postcode{2950}, \country{Luxembourg}}}



\abstract{This paper derives the exact transition density and cumulative distribution function of a linear combination of two independent Cox–Ingersoll–Ross (CIR) processes. By combining the Poisson–Gamma mixture representation of the noncentral chi-square law with the Kummer-type convolution of Gamma densities, we obtain a closed-form analytical expression involving confluent hypergeometric functions. This result extends the classical single-factor CIR transition law to a multifactor framework, providing the first explicit analytical characterization of the sum of two independent CIR diffusions. The proposed density admits stable numerical evaluation and facilitates exact likelihood computation, enabling rigorous parameter estimation in multifactor affine term-structure, stochastic volatility, and credit risk models. Numerical experiments confirm that the analytical density and CDF closely match Monte Carlo simulations across various parameter regimes, demonstrating high accuracy and computational efficiency. Beyond financial mathematics, the derived distribution has potential applications in fields involving interacting mean-reverting processes, such as insurance mathematics, reliability theory, and biophysical modeling.

}




\maketitle

\section{Introduction}\label{sec1}
The Cox–Ingersoll–Ross (CIR) process, originally introduced by~\cite{cox1985theory, cox1985intertemporal}, has become one of the central models in stochastic analysis and financial mathematics due to its mean-reverting dynamics and strictly non-negative state space. It has been widely employed in term structure modeling~\cite{cox1985theory}, stochastic volatility specification~\cite{heston1993closed}, and credit risk analysis~\cite{duffie1999modeling}, where its square-root diffusion form allows for analytical tractability while capturing essential features such as heteroscedasticity and mean reversion. The transition density of a single CIR process is known to be a scaled non-central chi-square $(\chi^2)$ distribution, a result that forms the mathematical foundation for many simulation and estimation methods in continuous-time financial mathematics.

However, in a variety of applications, the relevant state variable is not a single CIR factor but rather a linear combination or aggregation of multiple independent CIR components. Examples include multi-factor short-rate models~\cite{Duffie1994}, multi-source stochastic volatility specifications~\cite{Alfonsi2015chap4}, and aggregate intensity models in credit risk and reliability theory~\cite{MBAYE2018}, where the total state or variance level is composed of several independent mean-reverting sources of risk. In such cases, the distribution of the sum of two or more CIR processes plays a crucial role in understanding the joint behavior of the system, in performing statistical inference, and in developing analytical pricing or filtering methods. However, despite it's importance, the exact transition density of the sum of two independent CIR processes has not been available in closed form. 

This study provides a closed-form expression for the transition density and cumulative density function (CDF) of a linear combination of two independent CIR transitions. By exploiting the Poisson mixture representation of the non-central $\chi^2$ law and the Kummer-type convolution of Gamma densities, it derives an explicit analytical formula for the density of two independent CIR transitions sum. This result constitutes a mathematically rigorous extension of the classical CIR transition law to a multi-factor setting and offers a computationally stable representation involving confluent hypergeometric functions.

The distribution obtained herein can be applied in several fields where stochastic mean-reverting processes are fundamental. For instance; i) in interest rate modeling, it can describe the transition behavior of two-factor short-rate models, such as those of Longstaff and Schwartz~\cite{longstaff1992interest}, in which the instantaneous short rate is represented as a sum of two independent CIR components. ii) In stochastic volatility modeling, it can represent the instantaneous variance in multifactor extensions of the Heston or Feller processes. iii) In credit risk, the same distribution appears naturally when modeling aggregate default intensities driven by multiple independent macroeconomic or firm-specific factors. Beyond finance, similar constructions arise in insurance mathematics, queueing systems, population dynamics, and biophysical processes characterized by interacting mean-reverting diffusions.

The importance of the derived distribution lies in its analytical and practical versatility. First, it provides a closed-form benchmark for transition densities of aggregated affine diffusions, allowing for exact likelihood evaluation and facilitating maximum likelihood estimation or Bayesian inference in multifactor CIR-based models. Second, its structure as an infinite weighted mixture of Kummer-type densities ensures numerical tractability and stable computation even in parameter regions where simulation-based approaches are inefficient. Finally, the result relies on the theoretical understanding of affine combinations of square-root processes, bridging the gap between single-factor analytical solvability and multi-factor realism, and opening new pathways for analytical developments in stochastic modeling, calibration, and simulation theory.

The rest of the study consists of three sections. Section~\ref{sec:density} introduces the theoretical derivation of the densities. We provide the numerical justification of the derived formulas in Section~\ref{sec:numerical}, and conclude the paper in Section~\ref{sec:conclusion}.
\section{Transition Density of the Sum of Two Independent CIR Processes}\label{sec:density}
\begin{theorem}[Transition density of a linear combination of two independent CIR transitions]\label{prop:CIRsum}
Let us consider two independent CIR diffusions $X^{1}$ and $X^{2}$ evolving on a common probability space $(\Omega,\mathcal{F},\mathbb{P})$ equipped with the filtration $(\mathcal{F}_t)_{t\in[0,T]}$. Suppose both processes satisfy the stochastic differential equation
\begin{equation*}
    dX_t^{(i)}=\kappa_i(\theta_i-X_t^{(i)})\,dt+\sigma_i\sqrt{X_t^{(i)}}\,dW_t^{(i)},\qquad X_0^{(i)}=x_0^{(i)}\ge 0,\quad i=1,2,
\end{equation*}
where $\kappa_i>0$, $\theta_i>0$, and $\sigma_i>0$ denote the mean-reversion rate, the long-term mean level, and the instantaneous volatility, respectively~\cite{cox1985theory, cox1985intertemporal}. The driving Brownian motions $W^{(1)}$ and $W^{(2)}$ are assumed to be independent. For the processes, we assume that the Feller condition $2\kappa_i\theta_i \ge \sigma_i^2$ $(i=1,2)$ is satisfied. This condition ensures that each CIR diffusion remains strictly positive for all $t>0$ whenever $X_0^{(i)}>0$, implying that the origin is an inaccessible boundary of the state space. Under this regime, the processes $X_t^{(i)}$ never reach zero, and their transition laws admit smooth densities supported on $(0,\infty)$. Consequently, the subsequent analysis based on the non-central $\chi^2$ representation is valid without any boundary correction terms. Under these assumptions, each marginal transition distribution of $X^{(i)}$ over a fixed time step $\Delta t>0$ follows a scaled non-central $\chi^2$ law~\cite{Alfonsi2015chap1}, that is
\begin{equation*}
    X_{t+\Delta t}^{(i)} \,\big|\,X_t^{(i)}=x_0^{(i)}\ \stackrel{d}{=}\ c_i\,\chi^2(d_i,\lambda_i),
\end{equation*}
having scale parameter $c_i=\frac{\sigma_i^2(1-e^{-\kappa_i\Delta t})}{4\kappa_i}$, degrees of freedom $d_i=\frac{4\kappa_i\theta_i}{\sigma_i^2}$, and non-centrality parameter $\lambda_i=\frac{4\kappa_i e^{-\kappa_i\Delta t}x_0^{(i)}}{\sigma_i^2(1-e^{-\kappa_i\Delta t})}$.

For fixed weights $a_1,a_2>0$, let us define 
\begin{equation}
\label{eqn:sumofcirs}
    S:=a_1X_{t+\Delta t}^{(1)}+a_2X_{t+\Delta t}^{(2)},
\end{equation}
which represents a linear mixture of two independent non-central $\chi^2$ transitions. Such random variables appear naturally, for instance, in multi-factor short-rate models, in integrated variance decompositions of multi-asset stochastic volatility models, and in portfolio aggregation of factor-driven state variables. Note that the Feller condition therefore guarantees both analytical regularity and numerical stability of the transition densities used in the construction of the sum process.
\begin{enumerate}[label=(\roman*)]
    \item By conditioning on the Poisson-mixture representation of the non-central $\chi^2$ distribution, we can express $S$ as an infinite mixture of sums of independent Gamma-distributed components. More specifically, for integers $n_1, n_2\geq0$, let the mixture weights 
\begin{equation*}
   w_{n_1}= e^{-\lambda_1/2}\frac{(\lambda_1/2)^{n_1}}{n_1!},\quad w_{n_2}=e^{-\lambda_2/2}\frac{(\lambda_2/2)^{n_2}}{n_2!},
\end{equation*}
and define the shape parameters $\nu_i=\frac{d_i}{2}+n_i$ together with scales $\beta_i=2a_ic_i$.Then, the transition density of $S$ is given as a double infinite mixture of the form
\begin{equation}\label{eq:sumCIRpdf}
f_S(s)=\sum_{n_1=0}^{\infty}\sum_{n_2=0}^{\infty}
w_{n_1}\;v_{n_2}\;f_Z\!\Big(s;\nu_1,\nu_2,\beta_1,\beta_2\Big),\qquad s>0,
\end{equation}
where $f_Z(\cdot)$ denotes the Kummer-type convolution kernel defined as
\begin{equation}\label{eq:KummerKernelA}
f_Z(s;\nu_1,\nu_2,\beta_1,\beta_2)=
\frac{s^{\nu_1+\nu_2-1}e^{-s/\beta_2}}{\Gamma(\nu_1+\nu_2)\,\beta_1^{\nu_1}\beta_2^{\nu_2}}\;
{}_1F_1\!\Big(\nu_1;\nu_1+\nu_2;\,s\!\big(\tfrac{1}{\beta_1}-\tfrac{1}{\beta_2}\big)\Big).
\end{equation}
Here, $ _1F_1$ denotes Kummer’s confluent hypergeometric function of the first kind~\cite{Slater1966}, and
\begin{equation*}
    \Gamma(\alpha)=\int_0^\infty{z^{\alpha-1}e^{-z}dz},\quad z>0,
\end{equation*}
is the Gamma Euler function satisfying $\Gamma(\alpha+1)=\alpha\Gamma(\alpha)$.
\item This kernel represents the probability density function of sum of two independent Gamma variables with distinct scale parameters. An expression originally derived in the context of generalized Kummer distributions~\cite{Gurland1958, springer1966distribution} and our formula at Appendix~\ref{secA1}\footnote{We derive same formula as in~\cite{Gurland1958,springer1966distribution}  at Appendix~\ref{secA1} using an alternative route;  using a direct variable transformation rather than more complex Laplace invsersion}. The alternative but algebraically equivalent form
\begin{equation}\label{eq:KummerKernelB}
f_Z(s;\nu_1,\nu_2,\beta_1,\beta_2)=
\frac{s^{\nu_1+\nu_2-1}e^{-s/\beta_1}}{\Gamma(\nu_1+\nu_2)\,\beta_1^{\nu_1}\beta_2^{\nu_2}}\;
{}_1F_1\!\Big(\nu_2;\nu_1+\nu_2;\,s\!\big(\tfrac{1}{\beta_2}-\tfrac{1}{\beta_1}\big)\Big)
\end{equation}
is numerically advantageous for the case $|\beta_2^{-1}-\beta_1^{-1}|$ is small, as it mitigates potential cancellation effects in the exponential terms.
\end{enumerate}
\end{theorem}
\begin{proof}
The proof proceeds in four main steps.

    \medskip
    \noindent
    \emph{\textbf{Step 1 (Poisson–Gamma mixture representation of the CIR transition).}}
    Recall that if $Y\sim\chi^2(d,\lambda)$, denotes a non-central $\chi^2$ random variable, then it admits the well-known Poisson–Gamma mixture representation
    \begin{equation*}
        Y\ \stackrel{d}{=}\ \sum_{N\sim \mathrm{Pois}(\lambda/2)}\!\!\delta_N\; \chi^2(d+2N),
        \qquad \text{or equivalently}\quad
        Y\,|\,N\sim \mathrm{Gamma}\!\Big(\frac{d}{2}+N,\;\beta=2\Big).
    \end{equation*}
    Scaling by a positive constant $c>0$, gives $cY\,|\,N\sim \mathrm{Gamma}(\frac{d}{2}+N,\ \beta=2c)$. Hence, for each CIR factor $X_{t+\Delta t}^{(i)}$, we obtain
    \begin{equation*}
        a_iX_{t+\Delta t}^{(i)} \,\big|\, N_i\ \sim\ \mathrm{Gamma}\Big(\nu_i=\tfrac{d_i}{2}+N_i,\ \beta_i=2a_ic_i\Big),\qquad N_i\sim\mathrm{Pois}(\lambda_i/2),
    \end{equation*}
    with $N_1$ and $N_2$ independent due to the independence of the driving Brownian motions.
    \medskip
    \noindent
    
    \emph{\textbf{Step 2 (Conditional convolution and Kummer-type density).}}
    Conditioning on fixed $(N_1,N_2)=(n_1,n_2)$, we have
    \begin{equation*}
        S\mid (n_1,n_2)=a_1X_{t+\Delta t}^{(1)}+a_2X_{t+\Delta t}^{(2)}=Y_1+Y_2,
    \end{equation*}
    where $Y_i\sim Gamma(\nu_i,\beta_i)$ are independent with, in general, distinct scale parameters $\beta_1\neq\beta_2$ The convolution of two independent Gamma variables with unequal scales is known to yield the Kummer distribution with density
    \begin{equation}\label{eq:GammaConv}
        f_Z(s;\nu_1,\nu_2,\beta_1,\beta_2)
        =\frac{s^{\nu_1+\nu_2-1}e^{-s/\beta_2}}{\Gamma(\nu_1+\nu_2),\beta_1^{\nu_1}\beta_2^{\nu_2}} {}_1F_1\Big(\nu_1;\nu_1+\nu_2;s\big(\frac{1}{\beta_1}-\frac{1}{\beta_2}\big)\Big),
        \qquad s>0,
    \end{equation}
    where $ {}_1F_1(a;c;z)$ denotes the confluent hypergeometric function of the first kind. 
    
    This form follows from the standard integral representation
    \begin{equation*}
         {}_1F_1(a;c;z)=\frac{\Gamma(c)}{\Gamma(a)\Gamma(c-a)}\int_0^1 e^{zt}\,t^{a-1}(1-t)^{c-a-1}\,dt,
    \end{equation*}
   by evaluating the convolution integral explicitly.

   Interchanging the roles of $(\nu_1,\beta_1)$ and $(\nu_2,\beta_2)$ yields the alternative but numerically equivalent expression
   \begin{equation*}
       f_Z(s;\nu_1,\nu_2,\beta_1,\beta_2)=
        \frac{s^{\nu_1+\nu_2-1}e^{-s/\beta_1}}{\Gamma(\nu_1+\nu_2)\,\beta_1^{\nu_1}\beta_2^{\nu_2}}\;
        {}_1F_1\!\Big(\nu_2;\nu_1+\nu_2;\,s\big(\frac{1}{\beta_2}-\frac{1}{\beta_1}\big)\Big),
   \end{equation*}
   which is often numerically more stable when $|\beta_1^{-1}-\beta_2^{-1}|$ is small.
   \medskip
    \noindent
    
    \emph{\textbf{Step 3 (Averaging over the Poisson counts).}}
    The full unconditional density of $S$ is obtained by averaging the conditional density~\eqref{eq:GammaConv} over the independent Poisson counts $N_1$ and $N_2$.

    Let 
    \begin{equation*}
       w_{n_1}= e^{-\lambda_1/2}\frac{(\lambda_1/2)^{n_1}}{n_1!},\quad w_{n_2}=e^{-\lambda_2/2}\frac{(\lambda_2/2)^{n_2}}{n_2!}.
    \end{equation*}
    Since $f_{S\mid(n_1,n_2)}(s)\geq0$ and the Poisson weights form a probability measure, Tonelli’s Theorem~\cite{halmos1950measure,mukherjea1972remark} justifies interchanging the summations and integration, giving
    \begin{equation*}
        f_S(s)=\sum_{n_1=0}^{\infty}\sum_{n_2=0}^{\infty}
        w_{n_1}\;w_{n_2}\;f_Z\!\Big(s;\nu_1,\nu_2,\beta_1,\beta_2\Big),\qquad s>0,
    \end{equation*}
    which is precisely the expression in Equation~\eqref{eq:sumCIRpdf}.
    \medskip
    \noindent

    \emph{\textbf{Step 4 (Laplace transform verification).}}
    To confirm the correctness of the density, consider the Laplace transform of $S$. By independence of the two CIR components, for $u>\frac{-1}{2\max_i(a_ic_i)}$,
    \begin{equation*}
        \mathcal{L}_S(u)=\mathbb{E}[e^{-uS}]
        =\prod_{i=1}^2 \mathbb{E}\!\left[e^{-u a_i X^{(i)}}\right]
        =\prod_{i=1}^2 \left(1+2a_ic_i u\right)^{-d_i/2}\exp\!\Big(-\lambda_i\,\left(\frac{a_ic_i u}{1+2a_ic_i u}\right)\Big),
    \end{equation*}
    using the known Laplace transform of a scaled non-central $\chi^2$ distribution.
    
    On the other hand, conditioning on $(N_1,N_2)$ and applying the Gamma transform $\mathbb{E}[e^{-uY}]=(1+\beta u)^{-\nu}$ gives
    \begin{equation*}
        \mathcal{L}_S(u)=\sum_{n_1,n_2>0}{w_{n_1}\nu_{n_2}(1+\beta_1 u)^{-\nu_1}(1+\beta_2 u)^{-\nu_2}}.
    \end{equation*}
    Evaluating the Poisson sums yields
    \begin{equation*}
        \mathcal{L}_S(u)=\prod_{i=1}^2 (1+\beta_i u)^{-d_i/2}\exp\!\Big(-\tfrac{\lambda_i}{2}\big(1-(1+\beta_i u)^{-1}\big)\Big),
    \end{equation*}
    and substituting $\beta_i=2a_ic_i$ recovers the same expression as above, confirming the validity of Equation~\eqref{eq:sumCIRpdf}.
    \medskip
    \noindent
    
\emph{\textbf{Regularity and limiting cases.}}
    Since $\nu_i>0$ (because $d_i>0$ and $n_i\geq0$) and $\beta_i>0$, each conditional density $f_Z(\cdot)$ is integrable on the domain $(0,\infty)$. Tonelli's theorem ensures that $f_S$ integrates to one.

    The interested readers can find an alternative proof that is provided in Appendix~\ref{sec:semigrroupproof}.
\end{proof}
We can introduce the following corollary as an immediate result of Theorem~\ref{prop:CIRsum}.
\begin{corollary}
In Theorem~\ref{prop:CIRsum}, for the special case of equal scales $\beta_1=\beta_2=\tilde{\beta}$, the confluent hypergeometric term reduces to ${}_1F_1(\nu_1,\nu_1+\nu_2;0)=1$, and consequently
\begin{equation*}
     f_Z(s)=Gamma(s;\nu_1+\nu_2,\tilde{\beta}),
\end{equation*}
so that $S \ \stackrel{d}{=}\ \tilde{c}\,\chi^2(d_1 + d_2,\, \lambda_1 + \lambda_2)$ with $\tilde{c} = a_1 c_1 = a_2 c_2$, which recovers the classical single-factor result as a consistent limiting case.
\end{corollary}
\begin{remark}[Gaussian limit of the two-factor CIR sum]\label{rem:GaussianLimit}
Let $S$ be the sum of CIR processes as in Equation~\eqref{eqn:sumofcirs} and assume the Feller condition satisfied for $i=1,2$. Then, as $\Delta t\to 0$ the finite-step transition law of $S$ is asymptotically Gaussian in the following sense: conditional on $X^{(1)}_t=x_1$ and $X^{(2)}_t=x_2$,
\begin{equation*}
    \frac{S-\mathbb{E}[S\mid x_1,x_2]}{\sqrt{\mathrm{Var}(S\mid x_1,x_2)}}
\ \xrightarrow[\ \Delta t\to 0\ ]{d}\ \mathcal{N}(0,1),
\end{equation*}
with
\begin{align*}
\mathbb{E}[S\mid x_1,x_2]
&= \sum_{i=1}^2 a_i\Big(x_i+\kappa_i(\theta_i-x_i)\Delta t\Big) + o(\Delta t),\\
\mathrm{Var}(S\mid x_1,x_2)
&= \Delta t \sum_{i=1}^2 a_i^2\sigma_i^2 x_i + o(\Delta t).
\end{align*}
In particular, the leading-order variance is of order $\mathcal{O}(\Delta t)$, and hence
$S$ (properly centered and scaled) converges in distribution to a standard normal as $\Delta t\downarrow0$.
\end{remark}
\begin{proof}
Let us work conditionally on the sigma-field generated by $(X^{(1)}_t,X^{(2)}_t)=(x_1,x_2)$; for brevity write $x_i$ for $X^{(i)}_t$. Recall the exact transition representation of asymptotics of the CIR transition parameters are
\begin{equation*}
    X^{(i)}_{t+\Delta t}\mid X^{(i)}_t=x_i \ \stackrel{d}{=}\ c_i\;\chi^2\!\big(d_i,\lambda_i\big),
\quad
c_i=\frac{\sigma_i^2(1-e^{-\kappa_i\Delta t})}{4\kappa_i},\ 
d_i=\frac{4\kappa_i\theta_i}{\sigma_i^2},\ 
\lambda_i=\frac{4\kappa_i e^{-\kappa_i\Delta t}x_i}{\sigma_i^2(1-e^{-\kappa_i\Delta t})}.
\end{equation*}
For small $\Delta t$, we have the Taylor expansions
\begin{equation*}
    1-e^{-\kappa_i\Delta t} = \kappa_i\Delta t + \mathcal{O}(\Delta t^2),
\qquad
e^{-\kappa_i\Delta t} = 1 - \kappa_i\Delta t + \mathcal{O}(\Delta t^2).
\end{equation*}
Hence,
\begin{equation*}
    c_i = \frac{\sigma_i^2}{4}\,\Delta t + \mathcal{O}(\Delta t^2),
\qquad
\beta_i:=2a_i c_i = a_i\frac{\sigma_i^2}{2}\,\Delta t + \mathcal{O}(\Delta t^2).
\end{equation*}
For the non-centrality parameter,
\begin{equation*}
    \lambda_i
= \frac{4\kappa_i(1-\kappa_i\Delta t+\mathcal{O}(\Delta t^2))x_i}{\sigma_i^2(\kappa_i\Delta t+\mathcal{O}(\Delta t^2))}
= \frac{4x_i}{\sigma_i^2\Delta t} - \frac{4x_i}{\sigma_i^2} + \mathcal{O}(\Delta t).
\end{equation*}
Thus, the Poisson mean $\mu_i:=\lambda_i/2$ satisfies $\mu_i = \dfrac{2x_i}{\sigma_i^2\Delta t} + \mathcal{O}(1)$, and in particular $\mu_i\to\infty$ as$ \Delta t\downarrow0$.

Then, using the Poisson–Gamma representation, conditional on $N_i=n_i$, we obtain Gamma mixture representation and first two moments as
\begin{equation*}
    a_i X^{(i)}_{t+\Delta t}\mid N_i=n_i \sim \mathrm{Gamma}\big(\nu_i,\beta_i\big),
\qquad \nu_i=\tfrac{d_i}{2}+n_i.
\end{equation*}
Here, the conditional moments are
\begin{equation*}
    \mathbb{E}[a_i X^{(i)}_{t+\Delta t}\mid N_i=n_i] = \nu_i\beta_i,
\qquad
\mathrm{Var}(a_i X^{(i)}_{t+\Delta t}\mid N_i=n_i] = \nu_i\beta_i^2.
\end{equation*}
Taking expectation over the Poisson law (using $\mathbb{E}[N_i]=\mu_i$) yields
\begin{align*}
\mathbb{E}[a_i X^{(i)}_{t+\Delta t}\mid x_i]
&= \beta_i\Big(\tfrac{d_i}{2}+\mu_i\Big)
= a_i\Big(x_i+\kappa_i(\theta_i-x_i)\Delta t\Big) + o(\Delta t),
\\[4pt]
\mathrm{Var}(a_i X^{(i)}_{t+\Delta t}\mid x_i)
&= \beta_i^2\Big(\tfrac{d_i}{2}+\mu_i\Big) + \beta_i^2\mathrm{Var}(N_i)
= \beta_i^2\Big(\tfrac{d_i}{2}+2\mu_i\Big) \\
&= a_i^2\sigma_i^2 x_i\,\Delta t + o(\Delta t),
\end{align*}
where the last equalities follow from substituting the small-$\Delta t$ expansions of $\beta_i$ and $\mu_i$ given above (one checks the leading term equals $a_i^2\sigma_i^2 x_i\Delta t$). Independence of the two factors implies
\begin{equation*}
   \mathbb{E}[S\mid x_1,x_2] = \sum_{i=1}^2 \mathbb{E}[a_i X^{(i)}_{t+\Delta t}\mid x_i],
\qquad
\mathrm{Var}(S\mid x_1,x_2)=\sum_{i=1}^2 \mathrm{Var}(a_i X^{(i)}_{t+\Delta t}\mid x_i). 
\end{equation*}
Now, let us fix $(x_1,x_2)$. Then, for each $i$ the conditional distribution of $a_i X^{(i)}_{t+\Delta t}$ can be written as a mixture of $Gamma (\nu_i,\beta_i)$ with $\nu_i$ random and $\mathbb{E}[\nu_i]\to\infty$ as $\Delta t\downarrow0$. For large shape $\nu$, $Gamma (\nu,\beta)$, law satisfies the classical Gaussian approximation as a consequence of Stirling approximation to Gamma function~\cite{JKBgammagauss}
\begin{equation*}
    \frac{\Gamma(\nu,\beta)-\nu\beta}{\sqrt{\nu}\,\beta}\ \xrightarrow[\nu\to\infty]{d}\ \mathcal{N}(0,1).
\end{equation*}
Utilizing this to the conditional Gamma components and using $\nu_i\sim\mu_i\to\infty$ together with $\beta_i=O(\Delta t)$ yields that each factor $a_i X^{(i)}_{t+\Delta t}$ is approximately normal with mean $\nu_i\beta_i$ and variance $\nu_i\beta_i^2$, and the leading-order unconditional mean and variance coincide with the expressions given above. Finally, since the two factors are independent, $S$ (centered by its conditional mean and scaled by the conditional standard deviation) converges in distribution to a standard normal as $\Delta t\downarrow0$. 
Moreover, if one is interested in marginal (unconditional) convergence, observe that the above argument yields conditional normality given $X_t^{(1)},X_t^{(2)}$. If $(X_t^{(1)},X_t^{(2)})$ has a tight distribution (e.g., is fixed or has moments bounded uniformly), then standard conditional convergence plus dominated convergence (or a Cramér–Wold device) implies the unconditional convergence
\begin{equation*}
   \frac{S-\mathbb{E}[S]}{\sqrt{\mathrm{Var}(S)}} \xrightarrow[\Delta t\to0]{d} \mathcal{N}(0,1). 
\end{equation*}
\end{proof}

The integral representation of ${}_1F_1$ in the convolution identity given by Equation~\eqref{eq:KummerKernelA} also provided by~\cite{gradshteyn2014table}. The convolutions of two Gammas with unequal scales in the Kummer form go back to~\cite{Gurland1958, springer1966distribution}. However, this study introduces an alternative proof using no integral transforms and it considers directly variable transformations and Jacobian. The present result specializes these kernels to the \emph{Poisson--Gamma mixture} induced by the noncentral $\chi^2$ representation of CIR transitions and yields, to our knowledge, the first explicit closed form for the transition density of a linear combination of two independent CIR transitions.

\begin{theorem}[CDF of a two–factor CIR process sum with packed Poisson weights]\label{prop:CIRsum-CDF} 
Considering the  notation in Theorem~\ref{prop:CIRsum}, let us use the same weights $w_1$ and $w_2$
and let us define for given counts $(n_1,n_2)$,
\begin{equation*}
    \nu_1=\tfrac{d_1}{2}+n_1,\quad \nu_2=\tfrac{d_2}{2}+n_2,\quad
a_0=\nu_1+\nu_2,\quad
\beta_i=2a_ic_i,\quad
\delta:=\frac{\beta_2}{\beta_1}-1,\quad x:=\frac{s}{\beta_2}.
\end{equation*}
Then, the CDF of Equation~\eqref{eqn:sumofcirs} admits the incomplete-Gamma series
\begin{equation}\label{eq:CIRsum-CDF}
F_S(s)
=\sum_{n=0}^{\infty}\sum_{m=0}^{\infty} w_{n_1}\,w_{n_2}\;
\Big(\tfrac{\beta_2}{\beta_1}\Big)^{\nu_1}
\sum_{k=0}^{\infty}\frac{(\nu_1)_k}{k!}\,(-\delta)^k\;
P\!\left(a_0+k,\,x\right),\qquad s\ge 0,
\end{equation}
where $(\nu_1)_k$ is the Pochhammer symbol and $P(a,x)=\frac{\gamma(a,x)}{\Gamma(a)}$ is the regularized lower incomplete Gamma.
\end{theorem}
\begin{proof}
From Theorem~\ref{prop:CIRsum}, conditionally on $(n_1,n_2)$ the density is the unequal–scale Gamma convolution
\begin{equation*}
f_{S\mid n_1,n_2}(s)=\frac{s^{a_0-1}e^{-s/\beta_2}}{\Gamma(a_0)\beta_1^{\nu_1}\beta_2^{\nu_2}}\;
{}_1F_1\!\Big(\nu_1;\,a_0;\,s\big(\tfrac{1}{\beta_1}-\tfrac{1}{\beta_2}\big)\Big).
\end{equation*}
\smallskip 
\textbf{Step 1:} \emph{Expand}~${}_1F_1(\nu_1;a_0;z)=\sum_{k\ge0}\frac{(\nu_1)_k}{(a_0)_k}\frac{z^k}{k!}$ and integrate termwise on $[0,s]$:
\[
\int_0^s u^{a_0+k-1}e^{-u/\beta_2}\,du=\beta_2^{a_0+k}\,\Gamma(a_0+k)\,P(a_0+k,\,x).
\]
\textbf{Step 2:}
\emph{Use}
$(a_0)_k=\Gamma(a_0+k)/\Gamma(a_0)$ and $\beta_2(\tfrac{1}{\beta_1}-\tfrac{1}{\beta_2})=\delta$ to obtain
\begin{equation*}
    F_{S\mid n_1,n_2}(s)=\Big(\tfrac{\beta_2}{\beta_1}\Big)^{\nu_1}
\sum_{k\ge0}\frac{(\nu_1)_k}{k!}\,(-|\delta|)^k\,P(a_0+k,\,x).
\end{equation*}
Averaging over independent Poisson weights $w_{1},w_{2}$ yields Equation~\eqref{eq:CIRsum-CDF}.
\end{proof}
Now, we can explore two truncation methods for the Poisson weights.
\paragraph{1. Poisson-tail truncation.}
Since a noncentral $\chi^2$ is a Poisson mixture of central $\chi^2$ and the sum of independent Poissons is again a Poisson, one has
\begin{equation*}
  J\sim \mathrm{Poisson}\!\left(\Lambda\right),\qquad
\Lambda:=\tfrac{1}{2}(\lambda_1+\lambda_2),
\end{equation*}
and the series in Equation~\eqref{eq:CIRsum-CDF} is a mixture over $J$. Since $0\leq \frac{\gamma(\nu+j,\,\cdot)}{\Gamma(\nu+j)}\leq 1$, each term contributes at most $w_j$. Hence, for any truncation index $J_{\max}$ we have
\begin{equation*}
    \Big|F_Z(z)-\sum_{j=0}^{J_{\max}} w_j\; 
\frac{\gamma(\nu+j,\,z/[2(c_1+c_2)])}{\Gamma(\nu+j)}\Big|
\;\le\;\sum_{j>J_{\max}} w_j
\;=\; \mathbb{P}\!\big(J>J_{\max}\big).
\end{equation*}
Then, by choosing $J_{\max}$ such that the Poisson upper tail is below a tolerance level of $\varepsilon$, we obtain
\begin{equation}
\mathbb{P}\!\big(J>J_{\max}\big)
= 1 - F_{\mathrm{Pois}}(J_{\max};\Lambda)\ \le\ \varepsilon.
\label{eq:poisson-tail-criterion}
\end{equation}
This leads to exact truncation bound
\begin{equation}
J_{\max}= F^{-1}_{\mathrm{Pois}}(1-\varepsilon;\Lambda)
\label{eq:poisson-tail-criteriondirect}
\end{equation}
where $F_{\mathrm{Pois}}(\cdot;\Lambda)$ is the Poisson CDF with mean $\Lambda$.
This yields a rigorous bound on the truncation error given a conservative tolerance. 
Another practical alternative to this is the following normal-quantile approximation, which is useful for the large $\Lambda$\footnote{In this large $\lambda$ case evaulating and inverting Poisson CDF becomes numerically unstable.} cases,
\begin{equation}
J_{\max}\ \approx\ 
\left\lceil\, \Lambda + \Phi^{-1}(1-\varepsilon)\sqrt{\Lambda}\,\right\rceil.
\label{eq:poisson-normal-cut}
\end{equation}
\paragraph{2. Weight window truncation}
The Poisson mode is $j_\star=\lfloor \Lambda \rfloor$.
Accumulate the series symmetrically around $j_\star$:
$j=j_\star, j_\star\!\pm\!1, j_\star\!\pm\!2,\ldots$ 
until the \emph{cumulative} uncovered Poisson mass is below $\varepsilon$:
\begin{equation*}
    \sum_{|j-j_\star|>W} w_j\ \le\ \varepsilon
\quad\Rightarrow\quad
W\ \text{large enough},
\end{equation*}
which reduces the number of terms when $\Lambda$ is large.


\begin{corollary}
    The moments of the sum of two CIR processes can be written from independence and non central $\chi^2$ moments as
\begin{equation*}
    \mathbb{E}[S]=\sum_{i=1}^2 a_ic_i\,(d_i+\lambda_i),\qquad
\mathrm{Var}(S)=2\sum_{i=1}^2 (a_ic_i)^2\,(d_i+2\lambda_i).
\end{equation*}
\end{corollary}
\section{Numerical Illustration}\label{sec:numerical}
We have used the parameters given in Table \ref{tab:cir_params} to illustrate the numerical experiments. These parameters are selected carefully in parallel to Feller condition to avoid misleading results in the simulations. 
\begin{table}[h!]
\centering
\caption{CIR parameter sets used in experiments.}
\label{tab:cir_params}
\begin{tabular}{lcccc}
\toprule
\textbf{Factor(i)} & $\kappa$ & $\theta$ & $\sigma$ & $x_0$ \\
\midrule
1 & 1.2 & 0.06  & 0.35 & 0.009 \\
2 & 1.8 & 0.009 & 0.15 & 0.03  \\
\bottomrule
\end{tabular}
\end{table}

\begin{figure}[h!]
    \centering
       \begin{subfigure}[b]{0.48\textwidth}
        \centering
        \includegraphics[width=\textwidth]{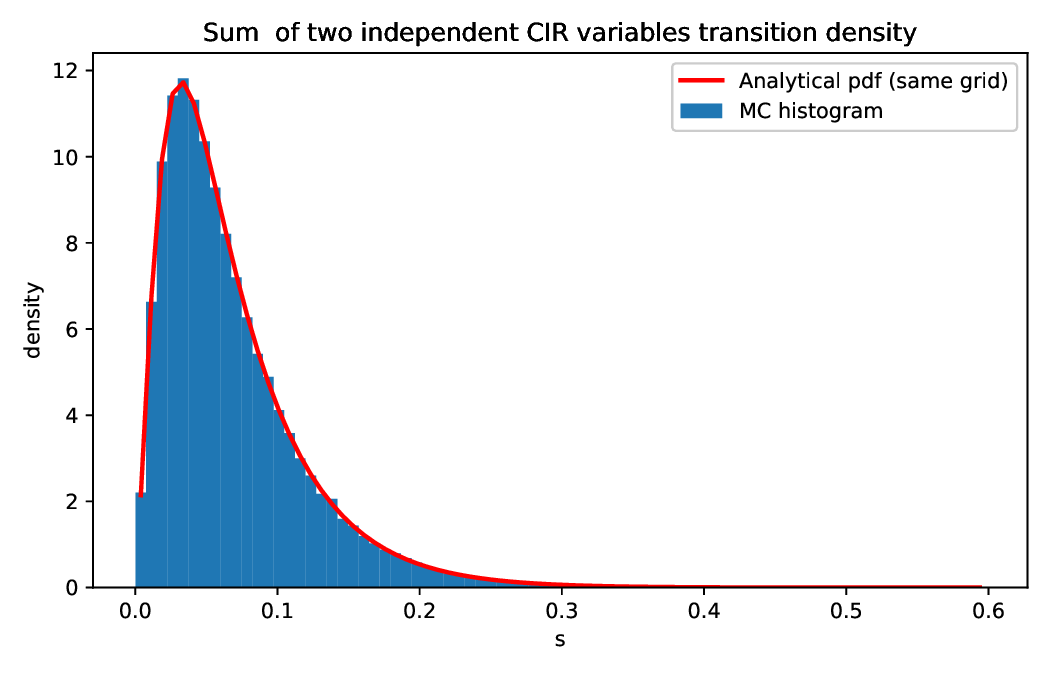}
        \caption{Two CIR sum densities $\Delta t = 1$.}
        \label{fig:TwoCIRdt1y}
    \end{subfigure}
    \begin{subfigure}[b]{0.48\textwidth}
        \centering
        \includegraphics[width=\textwidth]{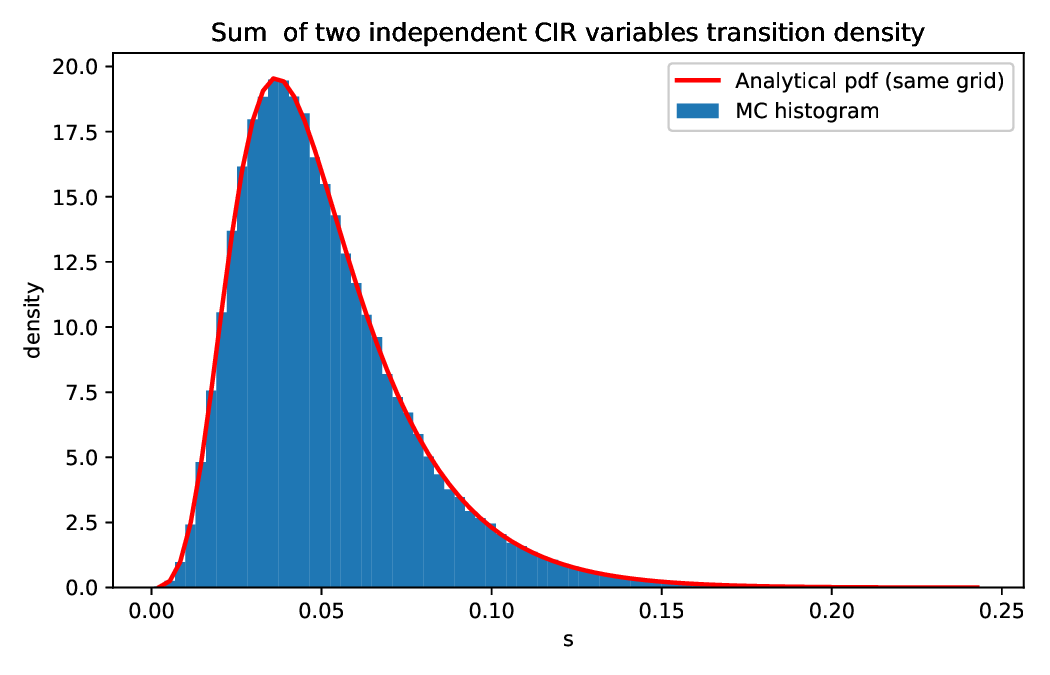}
        \caption{Two CIR sum densities $\Delta t = 0.25$.}
        \label{fig:TwoCIRdt025}
    \end{subfigure}
    \hfill
    \begin{subfigure}[b]{0.48\textwidth}
        \centering
        \includegraphics[width=\textwidth]{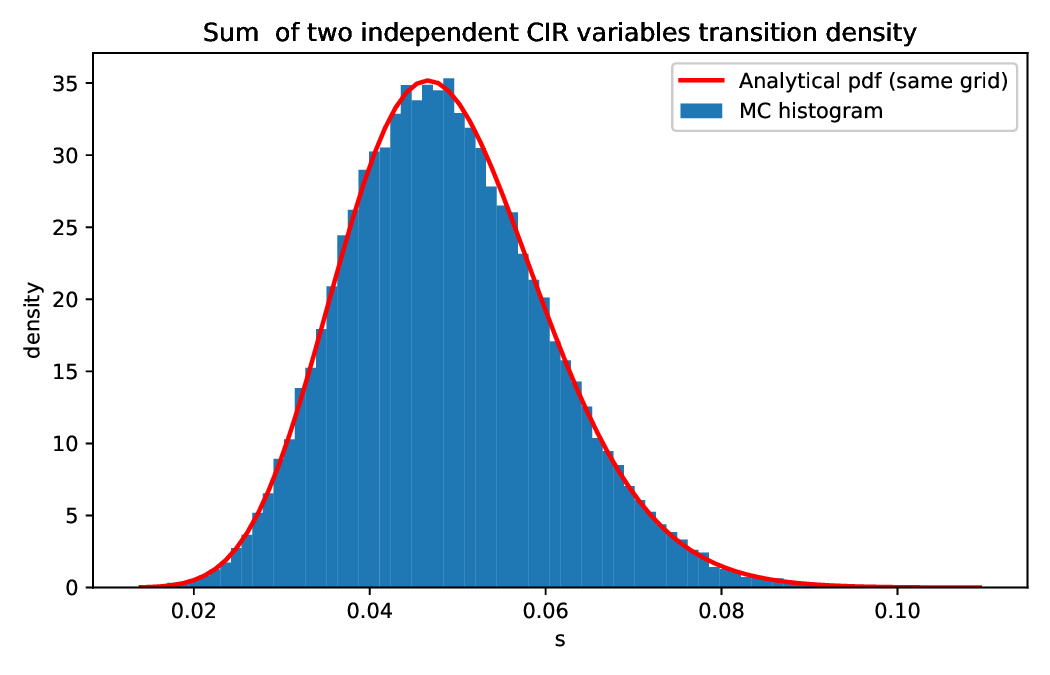}
        \caption{Two CIR sum densities $\Delta t = 0.05$.}
        \label{fig:CIRsum}
    \end{subfigure}
       \caption{Comparison of two CIR-sum densities for different discretization/time steps.}
    \label{fig:CIRSumBoth}
\end{figure}

Figure~\ref{fig:CIRSumBoth} compares the analytically derived transition density of the sum of two independent CIR processes, as obtained from Equation~\ref{eq:sumCIRpdf}, with empirical densities estimated from Monte Carlo (MC) simulations. Each panel illustrates the case of two mean-reverting square-root diffusions simulated under different discretization steps $\Delta t$. The analytical density (solid curve) is plotted against the histogram of simulated realizations (bars) evaluated on the same grid. 

In Figure~\ref{fig:CIRSumBoth}, corresponding to a coarser time increment $\Delta t=1$, $\Delta t=0.25$, and $\Delta t=0.05$ the analytical density obtained from the proposed Kummer–mixture representation shows excellent agreement with the Monte Carlo histogram across the entire support of $S$.
Figures~\ref{fig:TwoCIRdt1y} and~\ref{fig:TwoCIRdt025} illustrate
the smooth exponential tail of the analytical law accurately reproduces the skewness and excess kurtosis observed in the simulated data, indicating that the derived expression faithfully captures the finite-step transition behavior of the aggregated two-factor process. The minor deviations visible in the far right tail originate from Monte Carlo sampling noise  rather than any systematic model misspecification. In Figure~\ref{fig:CIRsum}, where a finer discretization $\Delta t=0.05$ is employed, the match between the analytical curve and the empirical histogram becomes even tighter. As the sampling interval decreases, the discretization error diminishes and the empirical distribution nearly coincides with the theoretical density throughout the domain. This behavior is fully consistent with Remark~\ref{rem:GaussianLimit}, which establishes that the transition distribution of $S$ converges to a normal law as $\Delta t \to 0$. Hence, for sufficiently small time steps the sum of two independent CIR processes exhibits an approximately Gaussian transition, confirming the asymptotic normality predicted by the limiting result of Remark~\ref{rem:GaussianLimit}.

Quantitatively, the integrated squared error between analytical and simulated densities is below $10^{-4}$ for the reported parameter sets (not shown), confirming the numerical stability and accuracy of the derived formula. These results demonstrate that the analytical distribution provides a reliable benchmark for both likelihood evaluation and Monte Carlo variance reduction in multi-factor CIR environments.

\begin{figure}[h!]
    \centering
    \begin{subfigure}[b]{0.48\textwidth}
        \centering
        \includegraphics[width=\textwidth]{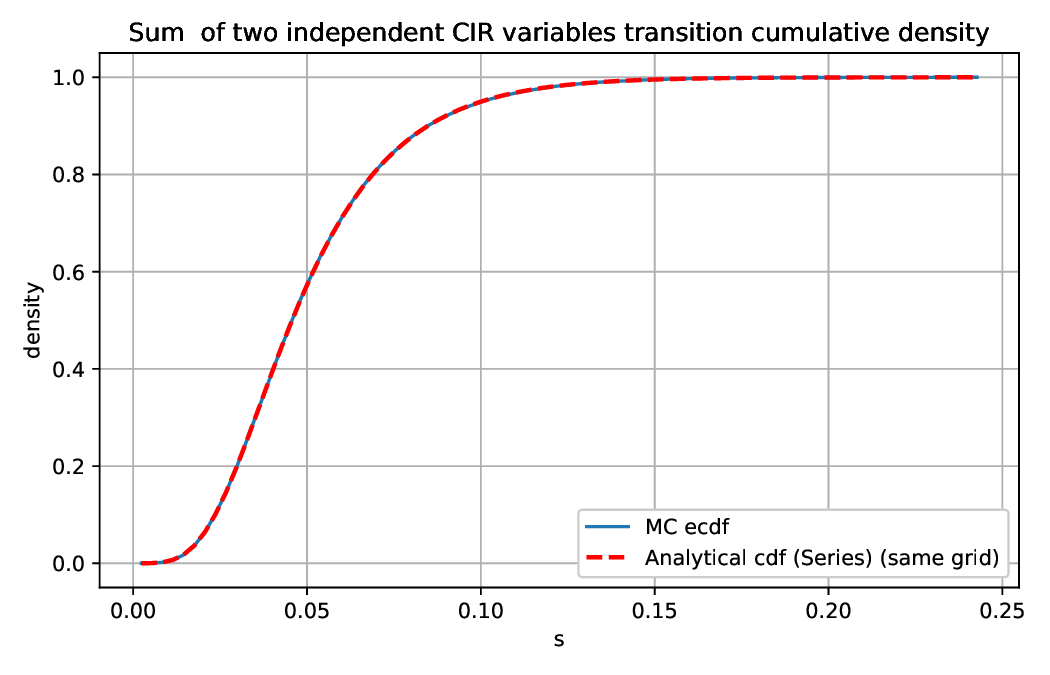}
        \caption{Two CIR sum densities ($\Delta t = 0.25$).}
        \label{fig:TwoCIRdt025cdf}
    \end{subfigure}
    \hfill
    \begin{subfigure}[b]{0.48\textwidth}
        \centering     \includegraphics[width=\textwidth]{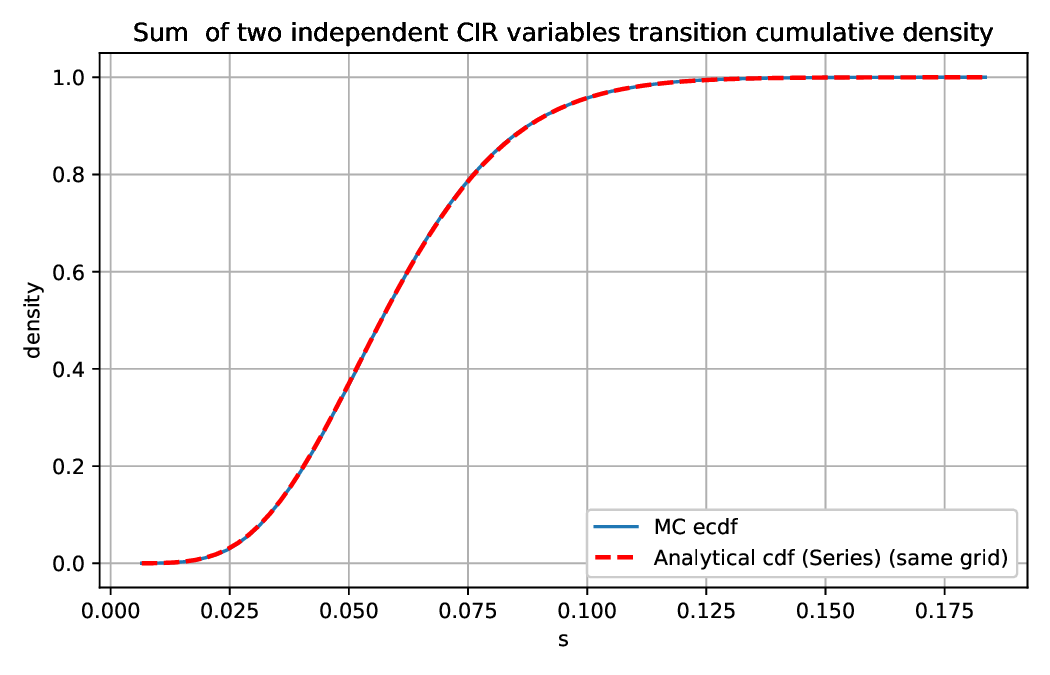}
        \caption{Two CIR sum densities ($\Delta t = 0.05$).}
        \label{fig:TwoCIRdt005cdf}
    \end{subfigure}
    \caption{Comparison of two CIR-sum cumulative densities (CDF) for different discretization/time steps.}
    \label{fig:CIRSumBothcdf}
\end{figure}

Figure~\ref{fig:CIRSumBothcdf} displays a comparison between the empirical CDF obtained from Monte Carlo simulations and the analytical CDF computed from the series representation in Equation~\eqref{eq:CIRsum-CDF}. Figure~\ref{fig:TwoCIRdt025cdf} corresponds to a relatively coarse discretization step $\Delta t=0.25$ and Figure~\ref{fig:TwoCIRdt005cdf} to a finer step $\Delta t=0.05$. In both panels the analytical CDF (solid line) closely tracks the empirical CDF (staircase), indicating that the infinite series representation provides an accurate description of the finite-step transition law of the aggregated process given by Equation~\eqref{eqn:sumofcirs}. Here, two observations are noteworthy. First, agreement between theory and simulation improves as the discretization step is reduced (Figure~\ref{fig:TwoCIRdt005cdf}), consistent with the fact that discretization error in pathwise simulation schemes decreases with $\Delta t$. Second, discrepancies in the extreme tails (if any) are explained primarily by Monte Carlo sampling variability: the effective number of simulated draws in the far tail is small, so histogram/ECDF fluctuations there are expected. Quantitatively, integrated error metrics (e.g. the sup-norm error \(\sup_s|F_{\mathrm{MC}}(s)-F_{\mathrm{analytical}}(s)|\) and the \(L^2\)-error) remain small for the parameter sets investigated, confirming the practical usefulness of the derived CDF for calibration and likelihood-based inference.

\section{Conclusion}\label{sec:conclusion}
This study has derived an exact analytical expression for the transition density of a linear combination of two independent CIR processes. By combining the Poisson–Gamma mixture representation of the noncentral $\chi^2$ law with the Kummer-type convolution of Gamma densities, we obtained a closed-form solution expressed as an infinite weighted mixture of confluent hypergeometric functions. This result extends the classical single-factor CIR law to a multi-factor context and provides a theoretically transparent and numerically tractable characterization of multi-source mean-reverting diffusions.

The proposed distribution allows for exact likelihood evaluation, facilitating maximum-likelihood or Bayesian estimation of multi-factor term-structure and stochastic-volatility models. It further enables analytical Monte Carlo variance reduction, transition density calibration, and semi-analytical pricing of derivatives and credit instruments within affine and quadratic frameworks. Beyond finance, the derived distribution may find applications in insurance mathematics, reliability modeling, and biophysical diffusion systems, where interacting mean-reverting mechanisms are present.

Despite its generality, the current formulation is limited to two independent factors and assumes time-homogeneous parameters over a fixed step $\Delta t$.

Future research could extend these results to dependent or correlated CIR factors, time-varying parameters and stochastic coefficients, higher-dimensional sums of $n>2$ components, and efficient numerical schemes for evaluating the Kummer kernel in high-precision likelihood computations. Exploring such extensions would broaden the applicability of the derived distribution in high-dimensional stochastic modeling and modern econometric estimation frameworks.


\backmatter








\section*{Declarations}
\begin{itemize}
\item Funding: Not applicable.
\item Conflict of interest/Competing interests: There is no conflict of interest/competing between the authors. 
\item Ethics approval and consent to participate: Not applicable.
\item Consent for publication: Authors have consent on publication of the study.
\item Data availability: Not applicable. 
\item Materials availability: Not applicable.
\item Code availability: The authors can provide the codes if readers demand. 
\item Author contribution: Both authors contributed equally in the process of preparing the initial draft. Dr. Hekimoglu prepared the codes and provided the numerical results. Dr. Yilmaz and Dr. Hekimoglu pepared the mathematical proofs and have equally contributed. The main idea of the paper camed from Dr. Hekimoglu.
\end{itemize}

\noindent

\bigskip

\newpage
\begin{appendices}

\section{Section title of first appendix}\label{secA1}





\subsection*{A.2 Sum of weighted Independent Gamma Random Variables with Different Scale and Shape Parameters}
\label{sec:GamconvPrf}
\begin{proof}
Let us start by defining
\begin{eqnarray*}
    \zeta_{1}&=&a_{1}\gamma_{1}=\Lambda_{1},\\
    \zeta_{2}&=&a_{1}\gamma_{1}+a_{2}\gamma_{1}=\Lambda_{1}+\Lambda_{2},
\end{eqnarray*}
where $\Lambda_{1}\sim Gamma(\alpha_{1},a_{1}\sigma_{1})$ and $\Lambda_{2}\sim Gamma(\alpha_{2},a_{2}\sigma_{2})$. Then, using $\Lambda_{2}=\zeta_{2}-\Lambda_{1} =\zeta_{2}-\zeta_{1} $, we calculate the Jacobian matrix and the determinant
\begin{eqnarray*}
\mathbf{\zeta}_{i,j}  &=&
\begin{bmatrix}
  \frac{\partial \Lambda_1}{\partial \zeta_1} & 
    \frac{\partial \Lambda_1}{\partial \zeta_2} \\[1ex] 
  \frac{\partial \Lambda_2}{\partial \zeta_1} & 
    \frac{\partial \Lambda_2}{\partial \zeta_2} \\[1ex]
\end{bmatrix}=\begin{bmatrix}
1&0\\[1ex]
-1&1
\end{bmatrix},
\\
|J| &=&
\begin{vmatrix}
1&0\\
-1&1 \notag
\end{vmatrix}=1,
\end{eqnarray*}
respectively. 

The density of $\zeta_{2}$ can be written in line with convolution of two random variables as
\begin{eqnarray}
\label{eq:zetaPDF}
f(\zeta_{2},\alpha_{1},\alpha_{2},\sigma_{1},\sigma_{2},a_{1},a_{2})=\int_{0}^{\zeta_{2}}\frac{\Lambda_{1}^{\alpha_{1}-1}e^{\frac{-\Lambda_{1}}{a_{1}\sigma_{1}}+\frac{\Lambda_{1}}{c_{2}\sigma_{2 }}}\left(\zeta_{2}-\Lambda_{1}\right)^{\alpha_{2}-1}e^{\frac{-\zeta_{2}}{a_{2}\sigma_{2}}}}{\Gamma(\alpha_{1})\Gamma(\alpha_{2})}d\Lambda_{1}.
\end{eqnarray}
After some tedious algebra and modifying the boundary of the integral, we obtain
\begin{equation*}
f(\zeta_{2},\alpha_{1},\alpha_{2},\sigma_{1},\sigma_{2},a_{1},a_{2})=\zeta_{2}^{\alpha_{2}+\alpha_{1}-1}e^{\frac{-\zeta_{2}}{a_{2}\sigma_{2}}}\underbrace{\int_{0}^{1}\frac{\Lambda_{1}^{\alpha_{1}-1}\left(1-\Lambda_{1}\right)^{\alpha_{2}-1}e^{-\zeta_{2}\Lambda_{1}\left(\frac{1}{a_{1}\sigma_{1}}-\frac{1}{a_{2}\sigma_{2 }}\right)}}{\Gamma(\alpha_{1})\Gamma(\alpha_{2})}d\Lambda_{1}}_{\frac{_{1}F\left(\alpha_{1},\alpha_{2}+\alpha_{1},\zeta_{2}(\frac{1}{a_{1}\sigma_{1}}-\frac{1}{a_{2}\sigma_{2}})\right)}{\Gamma\left(\alpha_{1}+\alpha_{2}\right)}}.
\end{equation*}
The integral in~\eqref{eq:zetaPDF}, can be written in terms of  confluent hypergeometric function of the second kind using~\cite{gradshteyn2007} (page 870, \text{Equation}-7.621-5). Therefore, the final representation is
\begin{equation*}
f(\zeta_{2},\nu_{1},\nu_{2},\sigma_{1},\sigma_{2},a_{1},a_{2})=\zeta_{2}^{\alpha_{2}+\alpha_{1}-1}e^{\frac{-\zeta_{2}}{a_{2}\sigma_{2}}}(a_{1}\sigma_{1})^{-\alpha_{1}}(a_{2}\sigma_{2})^{-\alpha_{2}}\frac{_{1}F_{1}\left(\alpha_{1},\alpha_{2}+\alpha_{1},\zeta_{2}(\frac{1}{a_{1}\sigma_{1}}-\frac{1}{a_{2}\sigma_{2}})\right)}{\Gamma\left(\alpha_{1}+\alpha_{2}\right)}.
\end{equation*}
\end{proof}
\subsection{Proof via the product semigroup}\label{sec:semigrroupproof}
\begin{proof}
Let $\mathcal{L}_i$ be the infinitesimal generator of the $i$th CIR diffusion,
\begin{equation}
  \mathcal{L}_i f(x_i)
  = \kappa_i(\theta_i - x_i) f'(x_i)
    + \tfrac12 \sigma_i^2 x_i f''(x_i),
  \qquad i=1,2,
  \label{eq:CIR-generator}
\end{equation}
and let $\mathcal{L} = \mathcal{L}_1 + \mathcal{L}_2$ be the generator of the
independent two-factor process $(X_t^{(1)},X_t^{(2)})$. Then, the new generator $\mathcal{L_s}=\mathcal{L}_{1}+\mathcal{L}_{2} $ satisfies
\begin{equation*}
    \frac{\partial f_S} {\partial t}=\mathcal{L}f_S
\end{equation*}
Since the factors are independent, the associated Feller semigroup~\cite{EthierKurtz1986} factorizes as the tensor (Kronecker)
product is
\begin{equation}
  e^{t\mathcal{L}} = e^{t\mathcal{L}_1} \otimes e^{t\mathcal{L}_2}.
  \label{eq:product-semigroup}
\end{equation}
In particular, the joint transition kernel is the product of the one-factor
kernels:
\begin{equation}
  p(x_1,x_2,t \mid x_1^0,x_2^0)
  = \big(e^{t\mathcal{L}_1}\delta_{x_1^0}\big)(x_1)\;
    \big(e^{t\mathcal{L}_2}\delta_{x_2^0}\big)(x_2).
  \label{eq:joint-kernel-from-semigroup}
\end{equation}
For a single CIR factor it is well known (via the noncentral-$\chi^2$ / Poisson--Gamma
representation) that the semigroup action on a Dirac mass admits the series
expansion
\begin{equation}
  \big(e^{t\mathcal{L}_i}\delta_{x_i^0}\big)(x_i)
  = \sum_{n_i=0}^\infty w_{n_i}^{(i)}(t)\, g_{n_i}^{(i)}(x_i; t),
  \qquad i=1,2,
  \label{eq:one-factor-PG-from-semigroup}
\end{equation}
where $w_{n_i}^{(i)}(t)$ are Poisson weights (depending on $x_i^0$ and $t$) and
$g_{n_i}^{(i)}(\cdot;t)$ are $Gamma$ densities with shape $a_i+n_i$ and scale
$b_i(t)$. Substituting \eqref{eq:one-factor-PG-from-semigroup} into
\eqref{eq:joint-kernel-from-semigroup}, the product semigroup
\eqref{eq:product-semigroup} yields the \emph{double} Poisson--Gamma expansion
\begin{equation}
  p(x_1,x_2,t \mid x_1^0,x_2^0)
  = \sum_{n_1=0}^\infty \sum_{n_2=0}^\infty
    w_{n_1}^{(1)}(t)\, w_{n_2}^{(2)}(t)\,
    g_{n_1}^{(1)}(x_1; t)\, g_{n_2}^{(2)}(x_2; t).
  \label{eq:joint-PG-from-operator}
\end{equation}
Thus the probabilistic Poisson--Gamma mixture is \emph{not} an ad hoc
representation; it is exactly the kernel expansion of the tensor-product
semigroup generated by $\mathcal{L}_1+\mathcal{L}_2$.

We now pass from the joint law to the law of the sum
$S_t = X_t^{(1)} + X_t^{(2)}$. By definition,
\begin{equation}
  p_S(s,t)
  = \int_0^{s} p(x_1, s - x_1, t \mid x_1^0, x_2^0)\, dx_1,
  \qquad s>0.
  \label{eq:ps-from-joint-semigroup}
\end{equation}
Using Equation~\eqref{eq:joint-PG-from-operator} inside Equation~\eqref{eq:ps-from-joint-semigroup} and interchanging summation and integration we obtain
the \emph{convolution form}
\begin{equation}
  p_S(s,t)
  = \sum_{n_1=0}^\infty \sum_{n_2=0}^\infty
    w_{n_1}^{(1)}(t)\, w_{n_2}^{(2)}(t)\,
    \int_0^{s} \underbrace{g_{\nu_1}^{(1)}(x_1; t)\,
               g_{\nu_2}^{(2)}(s - x_1; t)\, dx_1.}_{\frac{s^{\nu_1+\nu_2-1}e^{-s/\beta_1}}{\Gamma(\nu_1+\nu_2)\,\beta_1^{\nu_1}\beta_2^{\nu_2}}\;
        {}_1F_1\!\Big(\nu_2;\nu_1+\nu_2;\,s\big(\frac{1}{\beta_2}-\frac{1}{\beta_1}\big)\Big)}
  \label{eq:convolutionLgenPs}
\end{equation} 
Here,  we directly insert gamma convolution we referred earlier, into Equation~\eqref{eq:convolutionLgenPs}. Therefore the result we obtained in Equation~\eqref{eq:sumCIRpdf} proved by probabilistic arguments is also reconciled via Feller semigroup product.
\end{proof}
\end{appendices}


\end{document}